\documentclass[11pt]{amsart}
\usepackage{amsmath}
\usepackage{amsthm}
\usepackage{amssymb}
\usepackage{amsthm}
\usepackage{color}
\usepackage{pb-diagram}
\usepackage[all]{xy}
\def \C{{\mathbb C}}
\def \Z{{\mathbb Z}}
\def \K{{\mathcal K}} 
 
\def \D{{\mathbb D}} 
\def \N{{\mathbb N}} 
\def \R{{\mathbb R}}

\DeclareMathOperator{\ind}{ind}

\newtheorem{thm}{Theorem}[section]  
\theoremstyle{definition}
\newtheorem{dfn}[thm]{Definition}
\newtheorem{prop}[thm]{Proposition}

\newtheorem{lemma}[thm]{Lemma} 
\newtheorem{cor}[thm]{Corollary} 
 
\newtheorem{remark}[thm]{Remark}

\title{A secondary index for non-Fredholm operators associated with quantum walks} 
\author{ Toshikazu Natsume}
\address{Nagoya Institute of Technology, Nagoya, Japan}
\author{Ryszard Nest}
\address{University of Copenhagen, Copenhagen, Denmark}

\begin{document} 

\begin{abstract}
We study an analogue of chirality operators associated with quantum walks on the binary tree. For those operators we introduce a  K-theoretic invariant, an analogue of the index of Fredholm operators, and compute its values  in the ring of di-adic integers.
\end{abstract}
\footnote{This research is part of the EU Staff Exchange project 101086394 "Operator Algebras That One Can See". It was partially supported by the University of Warsaw Thematic Research Programme "Quantum Symmetries".}
\maketitle
\section{introduction}
The notion of quantum walk seems to have been explicitly introduced  first in 1993 by 
Y. Aharonov, L. Davidovich, and N. Zagury \cite{ADZ}. It originates from the idea of a quantum mechanical system on a lattice with internal degrees of freedom and has some rather unexpected properties.  An example is that while for a classical random walk the variance of the position is of the order of square root of time whereas , in the quantum case it is of the order of time. The original paper \cite{ADZ} gives a concrete application of the observation. Since then the notion of quantum walk has become important in quantum algorithms and quantum information theory (see {\it e.g.} \cite{K}) .

A self-adjoint operator $\Gamma$ on a Hilbert space is called a {\it symmetry} if $\Gamma^2=1$. Notice that $\frac{1\pm \Gamma}{2}$ are projections.
Suppose that a pair of symmetries $\Gamma_1, \Gamma_2$ on a Hilbert space $H$ is given. The operator $U=\Gamma_1\Gamma_2$ is a unitary, and 
 $\Gamma_2 U\Gamma_2=U^*$. Set $Q=U - U^*$. Then $\Gamma_2Q+Q\Gamma_2=0$. This means that with respect to the decomposition $H= \textup{Ran}(\frac{1+\Gamma_2}{2}) \oplus \textup{Ran}(\frac{1-\Gamma_2}{2})$ we have
 $$ Q = \begin{pmatrix}
 0 & -Q_+^* \\
 Q_+ & 0 \\
 \end{pmatrix}.$$
 The target of the investigation is the operator $Q_+$ for specific symmetries arising in the study of quantum walks.  

\subsection{One-dimensional case}Let us review classical one-dimensional quantum walks. The Hilbert space $\ell^2(\Z)$ is equipped with a canonical orthonormal basis $\{e_n\}_{n\in\Z} $. Let $L$ be the forward-shift $Le_j=e_{j+1}$. The first symmetry is the operator $\Gamma$ on $\ell^2(\Z)\otimes\C^2$ defined by
\begin{equation*}
\Gamma = \frac{1}{\sqrt{2}}\begin{pmatrix}
1 & L^* \\
L & -1 \\
\end{pmatrix}.
\end{equation*}
The two-sphere $S^2$ can be described as the set $\{ (x, \zeta ) \in \R\times \C : x^2+|\zeta|^2 = 1 \}$. We call a sequence $(a, b) = (a(n), b(n)) \in S^2$ a {\it walk}. A walk $(a, b)$ is {\it non-wandering} if the limits
$$a(\pm\infty)=\lim_{n\rightarrow \pm\infty}a(n), \;b(\pm\infty)=\lim_{n\rightarrow \pm\infty}b(n)$$
exist. The two sequences $a, b$ define diagonal operators on $\ell^2(\Z)$ in a canonical way. Then the non-wandering walk $(a,b)$ defines an operator $C$ called a {\it coin operator} by 
\begin{equation*}
C = \begin{pmatrix}
a & b^*\\
b & -a \\
\end{pmatrix}.
\end{equation*} The coin operator is the second symmetry we use. 
The operator
\begin{equation*}
Q_+ = \tfrac{1-C}{2}Q\tfrac{1+C}{2} : \textup{Ran}\bigl(\tfrac{1+C}{2}\bigr) \rightarrow \textup{Ran}\bigl(\tfrac{1-C}{2}\bigr)
\end{equation*}
is the {\it chirality operator} for a given pair of symmetries $\Gamma, C$. The Suzuki-Tanaka index theorem \cite{ST} is the folowing.
\begin{thm}
The chirality operator $Q_+$ is Fredholm if and only if $|a(\pm\infty)| \ne \frac{1}{\sqrt{2}}$, and the Fredholm index is given by
 \begin{equation*}
 \ind Q_+  = \left\{\begin{array}{cl}
1 & \mbox{if}\; |a(+\infty)| < \frac{1}{\sqrt{2}} <  |a(-\infty)| \\
-1 & \mbox{if}\; |a(-\infty)| < \frac{1}{\sqrt{2}} <  |a(+\infty)| \\
0 & \mbox{otherwise}
\end{array}\right..
 \end{equation*}
\end{thm}
Their proof is by explicit computations of the dimensions involved. The reference \cite{MN} provides a more general proof based on noncommutative geometry.

\subsection{Quantum walk on the binary tree}In this paper we study an analogue of the Suzuki-Tanaka index theorem, where the the lattice $\Z$ is replaced by the binary tree $T$ (Section 2). The space of infinite paths on $T$ provides a compactification  $\overline{T}$ of the space $v(T)$ of vertices of the tree with the boundary $K=\overline{T}\setminus T$, a Cantor set. The shift on $\Z$ is replaced by a rescaled left shift $L$ on $T$, the associated symmetry $\Gamma$ has the form
$$
\Gamma=\begin{pmatrix}
p&\overline{q}L^*\\
qL&1-p-LL^*
\end{pmatrix}
$$ 
and the coin operator $C$ has the same form as before, {\it i.e.} a walk is given by a function
$$
v(T)\rightarrow S^2
$$ 
A walk is non-wandering if $C$ extends to a continuos function on $\overline{T}$. The associated chirality operator $Q_+$  is, in this case, never Fredholm but the following holds.

\begin{thm}
When
$$
|a(x)|-|p|\neq 0 \mbox{ for all }x\in K,
$$ 
the index of $Q_+$ is well defined as an element  $\partial Q_+\in K_0(C(K))$ and, for any Radon measure $\tau$ on $K$, 
$$
<\partial Q_+,\tau>=\tau (\{x\in K ; a(x)>|p|\})-\tau(\{x\in K ; a(x)<|p|\}).
$$

\end{thm} 

The following two sections set the notation, the index theorem itself is proved, for $p=0$, in the section 4 while the general case is the content of the section 5.

\section{The binary tree}
Let $T$ be the standard binary tree, and denote by $v(T), e(T)$ the sets of verteces and edges, respectively. We regard $T$ as an oriented tree with the root $r$. If $e : u \mapsto v$ is an edge, $u$ is the {\it parent} of $v$, and $v$ is a {\it child} of $u$. We write $\{u\} = par(v)$. The set of children of $u$ is denoted $chi(u)$. If $v, v' \in chi(u)$, then they are {\it siblings}. 

We consider the Hilbert space $\ell^2(T):= \ell^2(v(T))$. 

\begin{dfn}[{\cite [Def. 5.1]{W}}]
The shift operator $S$ on $\ell^2(T)$ is defined by 
$$ (Sf)(v) = \left\{\begin{array}{cl}
f(par(v))\;, & v\neq r \\
0\; ,& v = r
\end{array}\right..$$
\end{dfn}

The shift operators defined for general trees are not necessarily bounded. 

\begin{thm}[{\cite[Thm.1.20]{W}}]
The shift operator $S$ on a rooted tree $T$ is bounded if and only if  there exists a constant $B > 0$ such that $\sharp(chi(u))\leq B$ for all $u\in v(T)$. 
\end{thm}
For the standard binary tree we have:

\begin{prop}
The adjoint of $S$ is
$$S^*f(u) = \sum_{v \in chi(u)}f(v).$$
\end{prop}
\begin{proof}
Straightforward computations.
\end{proof}
We have that
$$(S^*S)(f)(u) = \sum_{v\in chi(u)}(Sf)(v) = \sum_{v\in chi(u)}f(par(v)) = 2f(u).$$
So, if we set $L=\frac{1}{\sqrt{2}} S$, then $L$ is an isometry. The structure of isometries on Hilbert spaces is known.
\begin{thm}[Wold-von Neumann decomposition]
If $U$ is an isometry on a Hilbert space, then $U$ is a unitary, or a direct sum of copies of the unilateral shift, or a direct sum of a unitary and copies of the unilateral shift.
\end{thm}
A consequence of the Wold-von Neumann Theorem is that every non-unitary isometry has  spectrum the closed unit disk. The unitary part can be described as follows.

For an isometry $U$ on a Hilbert space, denote by $Z$ the orthogonal complement of $\cup_n \ker ((U^*)^n)$. Then $U_{|Z}$ is a unitary. When $Z=0$, the isometry is called {\it proper}. A proper isometry is unitarily equivalent to an isometry of the form $V\otimes I$, where $V$ is the unilateral  shift on $\ell^2(\Z_+)$. 

\begin{lemma} \label{def:L}
The isometry $L=\frac{1}{\sqrt{2}}S$ is proper.
\end{lemma}
\begin{proof}
Recall that  $L^* e_u = (1/\sqrt{2})e_{par(u)}, u\neq r$, and $L^*e_r=0$. Then for any $u \in v(T)$ there exists an $n\in\N$ such that $(L^*)^n e_u=0$. This implies that $Z=(\cup_n \ker ((L^*)^n))^\perp = 0$.
\end{proof}
\begin{cor}
The $C^*$-algebra $C^*(L)$ generated by $L$ on $\ell^2(T)$ is isomorphic to the $C^*$-algebra $C^*(V)$ generated by $V$ on $\ell^2(\Z_+)$.
\end{cor}

We want to construct an operator analogous to $\Gamma$ in the classical setting. Set $E=[L^*, L] = 1-LL^*$, and
$$\Gamma = \begin{pmatrix}
0 & L^* \\
L & E \\
\end{pmatrix}.$$
Then $\Gamma$ is a symmetry on $\ell^2(T)\otimes\C^2$. We want to construct the coin operator $C$.  It is known that there exists a compactification $\overline{T}$  of $T$  such that $\overline{T}\setminus T = K$ the ternary Cantor set. We consider the closed subset $T_0 =v(T)\cup K$ of $\overline{T}$. Let $C$ be as above with $a,b \in C(T_0)$. Recall that $f\in C(T_0)$ acts on $\ell^2(T)$ as point-wise multiplication operator.

 We follow the line of constructions of  $U=\Gamma C$ and  $Q= U-U^*$.
Set
$$\varepsilon =\frac{1}{\sqrt{2}} \begin{pmatrix}
\sqrt{1+a} & -\sqrt{1-a} \\
b/\sqrt{1+a} &b/\sqrt {1-a} \\
\end{pmatrix}. $$
Then $\varepsilon$ is a unitary, and $\varepsilon^*C\varepsilon = \begin{pmatrix}
1 & 0 \\
0 & -1 \\
\end{pmatrix}$.
\begin{lemma}
We have $EL=0$.
\end{lemma}
\begin{proof}
$((1-F)S\phi)(u) = (S\phi)(u)-(S\phi(u') = \phi(par(u))-\phi(par(u'))= 0$.
\end{proof}
Using Lemma 2.7 repeatedly we can show that
\begin{equation*}
\varepsilon^*Q\varepsilon =\begin{pmatrix}
0 & -Q_+^* \\
Q_+ & 0 \\
\end{pmatrix},
\end{equation*}
where 
$$Q_+ = \frac{\overline{b}}{\sqrt{1-a}}\,L\sqrt{1+a} - \sqrt{1-a}\,L^*\frac{b}{\sqrt{1+a}} +\frac{\overline{b}}{\sqrt{1-a}}E\frac{b}{\sqrt{1+a}}.$$
We want to determine if $Q_+$ is Fredholm.

\section{Noncommutative geometric approach}

Denote by $A$ the $C^*$-algebra generated by $S, C(T_0)$ on $\ell^2(T)$. Then
\begin{prop}
The $C^*$-algebra $A$ contains the ideal $\mathcal{K}(\ell^2(T))$.
\end{prop}
\begin{proof}
It is enough to show that for any $u, v \in v(T)$, the rank one operator $\theta_{u,v} : \C e_v \rightarrow \C e_u$ belongs to $A$. There exists a unique path $\{u_j\}$ such that $u_0= r , u_n= u$ and that $u_{j-1}=par(u_j)$. Then $S^ne_u = e_r$. Denote by $\chi_u \in C(\overline{T})$) the characteristic function for the set $\{u\}$. Then $\theta_{r,u}= S^n \chi_u$. Therefore $\theta_{r, u}\in A$. Obviously $\theta_{u, r} = \theta_{r,u}^* \in A$. Consequently $\theta_{u,v} = \theta_{u,r}\theta_{r,v} \in A$.
\end{proof}

Consider the short exact sequence:
$$0 \rightarrow \mathcal{K}(H) \rightarrow A \xrightarrow{\sigma} A/\mathcal{K}(H) \rightarrow 0.$$
Notice that the quotient algebra $ A/\mathcal{K}(H)$ is not abelian. 

The following proposition is the key.

\begin{prop}
We have that $A/\mathcal{K} \cong C^*(L)\otimes C(K).$
\end{prop}
\begin{proof}
We want to construct a linear map $s : C^*(L)\otimes C(K) \rightarrow A$ such that $\sigma\circ s $ is bijective.  
The set $K_0$ of endpoints of the removed open sets in the construction is dense in $K$. Therefore continuous functions are determined by the restrictions on $K_0$. 

For any $x\in K$ there exists a unique sequence of vertexes $\{u_j\}$ such that $u_0(x)=r, u_j (x)= \textup{par}(u_{j+1}(x))$, and $\lim_{j\rightarrow \infty}u_j(x) = x$ in $\overline{T}$. We use this picture to construct $s$.

Let $f\in C^*(L)\otimes C(K) = C(K, C^*(L))$. We proceed by induction to create $\tilde{f} \in A$ . The first step is to take care of two points $0,1$. Set $ \tilde{f}(u_j(0))= f(0) \in C^*(L)$, and $ \tilde{f}(u_j(1))= f(1) \in C^*(L)$. The second step is to take care of the points $1/3, 2/3$. We have the sequence $\{u_j(1/3)\}$. Now we know that $u_1(0) = u_1(1/3)$. So we set $ C^*(V)ilde{f}(u_j(1/3))= f(1/3), j\geq 2$, and $\tilde{f}(u_j(1/3))= f(0), j=0, 1$. Similarly we set $ \tilde{f}(u_j(2/3))= f(2/3), j\geq 2$, and $ \tilde{f}(u_j(2/3))= f(1), j=0, 1$. The next step is to take care of $1/9, 2/9, 7/9, 8/9$. Then $u_2(1/9) = u_2(0)$. So set $ \tilde{f}(u_j(1/9))= f(1/9), j\geq 3$, and $ \tilde{f}(u_j(1/9))= f(0), j\leq 2$.

 We proceed by induction to define $ \tilde{f}(u)$. Notice that any $u \in v(T)$ is on a path to some $x \in K_0$ .  Consequently we can define $\tilde{f}\in A$.
 
 By the construction above it is obvious that the map $s(f) = \tilde{f}$ is linear, and $s(f)s(g) - s(fg) \in \K(\ell^2(T))$. So $\sigma\circ s$ is a $\ast$-homomorphism $C(K)\otimes C^*(L) \rightarrow A/\K(\ell^2(T))$. We need to show $\sigma\circ s$ is bijective. Let $f\cdot w , f\in C(T_0), w\in C^*(L)$. Then $s(f_{|K}) -f \in \K(\ell^2(T))$. This implies that $\sigma\circ s$ is surjective.  It is easy to see that if $s(f) \in \K(\ell^2(T))$, then $f=s(f)_{|K} = 0$.
\end{proof}

By Corollary 2.6, 
$$A/\mathcal{K}(\ell^2(T)) \cong C^*(V)\otimes C(K).$$

\noindent
{\it Observation.}  The operatot $Q_+$ is not Fredholm. 

In order to see this, we only have to show that $\sigma(Q_+)$ is not invertible in  $C^*(V)\otimes C(K)$. We have
$$\sigma(Q_+) = \frac{\sqrt{1+a}}{\sqrt{1-a}}\overline{b}V - \frac{\sqrt{1-a}}{\sqrt{1+a}}bV^* + |b|(1-VV^*).$$

We claim that, for each fixed values of $a(x), b(x), x\in K$, as an operator on $\ell^2(\Z_+)$, the operator $\sigma(Q_+)$ is not invertible. For this we show that $\ker(\sigma(Q_+)^* )$ is non-trivial. 

Let $\xi=\sum_{n=0}^\infty C_n e_n  \in \ell^2(\Z_+)$. Suppose that $a(x) > 0$, then it is straightforward to see that $\xi\in \ker(\sigma(Q_+)^* )$ if and only if
\begin{align*}
C_1 & = -\frac{\sqrt{1-a(x)}}{\sqrt{1+a(x)}}\frac{|b(x)|}{b(x)} C_0 \\
C_{n+1}& = -\frac{1-a(x)}{1+a(x)}\frac{|b(x)|}{b(x)} C_{n-1}, n\geq 1.
\end{align*}
Notice that $0 < \frac{1-a(x)}{1+a(x)} < 1$. It follows that $\dim \ker(\sigma(Q_+)^* ) =1$.
By a similar argument, if $a(x) < 0$ then $\dim \ker(\sigma(Q_+ )) =1$

What we have seen so far is that we cannot extract any numerical invariant as long as we stick to the classical Fredholm index. We still want to get numerical information from the operatoe $Q_+$.
Thus we need to extend the notion of indices.

Denote by $[A, A]$ the commutator ideal of $A$. For nonzero $f \in C(T_0)$, the commutator $[S, f]$ is nonzero, and belongs to $\mathcal{K}(H)$. This means that $[A,A]\cap \mathcal{K}(\ell^2(T)) \neq 0$. Hence $[A,A]\cap \mathcal{K}(\ell^2(T)) $ is a non-trivial ideal of $\mathcal{K}(\ell^2(T)) $. Since $ \mathcal{K}(\ell^2(T)) $ is simple, we have $[A,A]\cap \mathcal{K}(\ell^2(T)) = \mathcal{K}(\ell^2(T))$. Hence $[A, A] \supset \mathcal{K}(\ell^2(T))$.

Now we have the following commutative diagram:
$$\begin{array}{ccccccccc}
0 & \rightarrow & \mathcal{K} & \rightarrow& A & \rightarrow &A/\mathcal{K} & \rightarrow & 0 \\
&& \cap && \parallel && \downarrow && \\
0 & \rightarrow & [A,A]& \rightarrow &A & \rightarrow &A/[A,A] & \rightarrow & 0 \\
\end{array},$$
where $\K=\K(\ell^2(T))$.

We know 
$$A/[A,A] \cong C^*(L)/[C^*(L), C^*(L)]\otimes C(K),$$
and $C^*(L)/[C^*(L), C^*(L)]\cong C^*(V)/[C^*(V), C^*(V)] \cong C(S^1)$. Therefore we get the exact sequences:
$$\begin{array}{ccccccccc}
0 & \rightarrow & \mathcal{K} & \rightarrow& A & \rightarrow &C^*(V)\otimes C(K) & \rightarrow & 0 \\
&& \cap && \parallel && \downarrow && \\
0 & \rightarrow & [A,A]& \rightarrow &A & \rightarrow &C(S^1)\otimes C(K) & \rightarrow & 0 \\
\end{array}.$$

Since $K_1(\K) = K_1(C^*(V)\otimes C(K))=0$, a part of the six-term exact sequence is:
\begin{equation}
0 \rightarrow K_0(\K) \rightarrow K_0(A) \rightarrow K_0(C^*(V)\otimes C(K)) \rightarrow 0.
\end{equation}
Also we get $K_1(A)= 0$.
Since $K_0(C^*(V)) = \Z$ is generated by the class of unit, $K_0(C^*(V)\otimes C(K)) \cong K_0(C(K))$. A well-known fact is that
$$ K_0(C(K)) \cong \varinjlim \Z^{2^n}.$$
This implies 
$$K_0(A) \cong \Z\oplus\varinjlim \Z^{2^n}.$$
Apply the six-term exact sequence to the second short exact sequence we get
$$\begin{array}{ccccc}
K_0([A,A])  & \rightarrow& \Z\oplus\varinjlim \Z^{2^n} & \rightarrow & \varinjlim \Z^{2^n} \\
\uparrow &&&& \downarrow  \\
 \varinjlim \Z^{2^n}& \leftarrow &0 & \leftarrow & K_1([A, A])  \\
\end{array}.$$
Call the quoient map $\hat{\sigma} : A \rightarrow A/[A, A] = C(S^1)\otimes C(K)$ {\it the secondary symbol map}. Recall that the connecting map $\delta_1 : K_1(A/[A, A])) \rightarrow K_0([A, A]) $ is called {\it the index map}. For an $a \in A$ if the element $\hat{\sigma}(a)$ is invertible, then its analytic index $\ind a \in K_0([A, A]) \cong \Z\oplus\varinjlim \Z^{2^n}$ is defined to be $\delta_1([\hat{\sigma}(a)])) \in \varinjlim \Z^{2^n} \subset K_0([A, A])$. There exists a surjective homomorphism  from the group $\varinjlim \Z^{2^n}$  to the following group $G$ of di-adic integers:
$$ G = \Bigl\{ \tfrac{m}{2^n}  : m \in \Z, n = 0, 1, 2,\cdots \Bigr \}.$$

Let us go back to $Q_+$. Denote by $z$ the canonical generator of $C(S^1)$, and $a_K, b_K$ the restrictions of $a, b$ onto $K$, respectively. Then
$$\hat{\sigma}(Q_+) =  \frac{\overline{b_K}}{\sqrt{1-a_K}}\,z\sqrt{1+a_K} - \sqrt{1-a_K}\,\overline{z}\frac{b_K}{\sqrt{1+a_K}} ,$$
because $E =1-LL^* \in [A, A]$.

\begin{lemma}
The symbol $\hat{\sigma}(Q_+)$ is invertible in $C(S^1)\otimes C(K)$ if and only if $a_K, b_K$ are invertible in $C(K)$.
\end{lemma}
\begin{proof}
If there exists $x\in K$ such that $b_K(x)=0$, then for all $z\in S^1$ we have $\hat{\sigma}(Q_+)(x,z) = 0.$ Thus $\hat{\sigma}(Q_+)$ is not invertible.

Now suppose that there exists an $x\in K$ such that $a_K(x)=0$. From this it follows $|b_K(x)|=1$.   We have $\hat{\sigma}(Q_+) = \overline{b_K(x)}z - b_K(x)\overline{z}$. As $|b_K(x)|=1$, we have $\hat{\sigma}(Q_+)(x, b_K(x)) = 0$. Thus $\hat{\sigma}(Q_+)$ is not invertible.

Conversely, suppose that $a_K(x) \neq 0, b_K(x)\neq 0 $ for all $x\in K$. Then
$$\Bigl\vert \frac{\sqrt{1+a_K(x)}}{\sqrt{1-a_K(x)}}\overline{b_K(x)}\Bigr\vert \neq \Bigl\vert\frac{\sqrt{1-a_K(x)}}{\sqrt{1+a_K(x)}}b_K(x)\Bigr\vert.$$
This implies that $\hat{\sigma}(Q_+)(x, z)\neq 0 $ for all $z\in S^1, x\in K$.
\end{proof}

\begin{dfn}
The class $\delta_1([\hat{\sigma}(Q_+)]) \in \varinjlim \Z^{2^n}\subset K_0([A, A])$ is called the {\it secondary index } of the chirality operator $Q_+$.
\end{dfn}

Our goal is to extract a numerical invariant from $\delta_1([\hat{\sigma}(Q_+)])\in \varinjlim \Z^{2^n}$ . 

\section{The index theorem}

We need a measure on the Cantor set in order to construct a trace on $C(K)$. Recall that $K$ is identified with the infinite product $\prod_{i=1}^\infty X_i$ with $X_i= \{0, 1 \}$ the two point space. Let $\mu_0$ be the discrete measure on $\{0, 1\}$ given by $\mu_0(\{0\})=\mu_0(\{1\})=\tfrac{1}{2}$, and $\mu$ be the product measure on $K$. Then $\mu$ is a probability measure on $K$ whose cumulative distribution function is the Cantor function (Devil's staircase). Denote by $\tau$ the trace on $C(K)$ defined by $\mu$. 

 Let $\varepsilon$ be the canonical densely defined cyclic 1-cocycle on $C(S^1)$. We consider the cup product $\varepsilon\sharp\tau$ on $C(S^1)\otimes C(K)$.

\begin{thm}
 There exists a $\ast$-homomorphism $\pi : [A, A] \rightarrow \K(\ell^2(\Z_+))\otimes C(K)$ such that the sequence
 $$0 \rightarrow \K(H) \rightarrow [A, A] \rightarrow \K(\ell^2(\Z_+))\otimes C(K)\rightarrow 0 $$
 is exact.
 \end{thm}
 \begin{proof}
 Let $B$ be the $C^*$-subalgebra of $A$  generated by $L, L^*$. Then for any $F \in C(T_0)$ and $\gamma \in B$ the commutator $[\gamma, F]$ belongs to $\mathcal{K}(H)$. So, what we need to analyze is the commutator ideal $[B,B]$ of $B$. It is straightforward to see that $[B,B]$ is generated by $E= [L^*, L] = 1-LL^*$. The commutator ideal $[C^*(V), C^*(V)]$ of $C^*(V)$ is exactly $\mathcal{K}(\ell^2(\Z_+))$. Then the conclusion follows from the diagram:
 $$\begin{array}{ccccccccc}
0 & \rightarrow & \mathcal{K} & \rightarrow& [A, A]  & \rightarrow &[C^*(V), C^*(V)] \otimes C(K) & \rightarrow & 0 \\
&& \parallel  & & \cap &&\cap && \\
0 & \rightarrow & \K & \rightarrow &A & \rightarrow &C^*(V)\otimes C(K) & \rightarrow & 0 \\
\end{array}.$$
 \end{proof}
 
 Set $\omega=\pi^*(\textup{Tr}\sharp \tau)$, where $\textup{Tr}$ is the canonical trace on $\K(\ell^2(\Z_+))$. Let $\varphi : K_0([A,A]) \rightarrow \R$ be the map obtained by the coupling with $\omega$, and let $\psi : K_1(C(S^1)\otimes C(K))) \rightarrow \R$ be the map obtained by the coupling with $\varepsilon \sharp \tau$.
 
 \begin{thm}[The Index Theorem] The following diagram is commutative:
 \begin{equation*}
\xymatrix{
K_1(C(S^1)\otimes C(K)) \ar[rr]^-{\delta_1} \ar[dr]_{\psi} & \ar@{}[d]|{\circlearrowright}& K_0([A,A]) \ar[dl]^-{\varphi} \\
& \R & }.
\end{equation*} 
\end{thm}

\begin{dfn}
The {\it numerical secondary index} of the chirality operator is 
$$\textup{s-ind} Q_+ = \varphi(\delta_1(\hat{\sigma}(Q_+)) \in \R.$$
\end{dfn}
 
Recall that since $\hat{\sigma}(Q_+)$ is invertible, $a_K$ is nonzero. Set $U_+ = \{ x \in K : a_K(x) > 0 \}, U_- = \{ x \in K : a_K(x) < 0 \}$. Then $K= U_+ \cup U_-$ a disjoint union. Denote by $\chi_+, \chi_- $ the characteristic functions of $U_+, U_-$, respectively.  The functions $\chi_\pm$ are continuous on $K$.
 
 \begin{cor}
We have that $\textup{s-ind} \,Q _+= \tau(\chi_+) - \tau(\chi_-)$.
\end{cor}

For the proofs of Theorem 4.2 and Corollary 4.4 we need some facts about the $K$-theory of the Cantor set. The n-th step of the construction of the Cantor set removes $2^{n-1}$ open intervals $(x_i, y_i), i=1, \cdots, 2^{n-1}, x_1<y_1< \cdots, x_{2^{n-1}} < y_{2^{n-1}}$. Set $y_0=0, x_{2^{n-1}+1}=1$. Denote by $f_j^{(n)}$ the characteristic function of $K\cap [y_j, x_{j+1}]$. Then the $K_0$-classes of $f_j^{(n)}$'s generate $K_0(C(K))$. As above $z$ is the canonical generator of $C(S^1)$. Then the classes $[z]\otimes[f_j^{(n)}] \in K_1(C(S^1))\otimes K_0(C(K))$ generate $K_1(C(S^1)\otimes C(K))$. This class is represented by the unitary $u = z f_j^{(n)} + (1-f_j^{(n)}) \in C(S^1)\otimes C(K)$.

We need to compute $\delta_1([u])\in K_0([A, A])$.
We review a very useful formula to describe the index map in the six-term exact sequence of K-groups.

Let $0 \rightarrow I \rightarrow A \xrightarrow{\pi} A/I \rightarrow 0 $ be a short exact sequence of $C^*$-algebras with $A$ being unital. We give a description of the index map $\delta_0 : K_1(A/I) \rightarrow K_0(I)$.

For an invertible $u \in A/I$, let $M, N\in A$ be such that $\pi(M) =u, \pi(N)= u^{-1}$. Then
$$ \Omega = \begin{pmatrix}
1-MN & M \\
2N-NMN & -(1-NM) \\
\end{pmatrix}$$
satisfies $\Omega^2 = 1$. Hence $\Omega$ is invertible. Now $g= \Omega\bigl(\begin{smallmatrix}
0 & 1 \\
1 & 0\\
\end{smallmatrix}\bigr)$
is invertible in $M_2(A)$ and $\pi(g) = \bigl(\begin{smallmatrix} u & 0 \\
0 & u^{-1}\end{smallmatrix}\bigr)$. We have
$$g = \begin{pmatrix}
M & 1-MN \\
-(1-NM) & 2N-NMN \\
\end{pmatrix}.$$
Following the standard recipe 
$$\delta_0([u]) = \bigl[g\bigl(\begin{smallmatrix}
 1 & 0 \\
0 & 0\\
\end{smallmatrix}\bigr)g^{-1}\bigr] - \bigl[\bigl(\begin{smallmatrix}
 1 & 0 \\
0 & 0\\
\end{smallmatrix}\bigr)\bigr].$$
Recall $e = g \bigl(\begin{smallmatrix}
 1 & 0 \\
0 & 0\\
\end{smallmatrix}\bigr)g^{-1} \in M_2(I^\sim)$, and $\hat{e} = e - \bigl(\begin{smallmatrix}
 1 & 0 \\
0 & 0\\
\end{smallmatrix}\bigr) \in M_2(I)$. By straightforward computations we get
$$\hat{e} = \begin{pmatrix}
- (1-MN)^2 & \ast\ast \\
\ast & (1-NM)^2 \\
\end{pmatrix}, $$
where $\ast, \ast\ast$ are nontrivial terms. We are interested in diagonal entries to compute a pairing with a cyclic cocycle.

In general we look for a cyclic even-cocycle $\omega$ on $I$, and we extract a numerical invariant:
$$\langle \omega, [e]- \bigl[ \bigl(\begin{smallmatrix}
 1 & 0 \\
0 & 0\\
\end{smallmatrix}\bigr)\bigr] \rangle = \omega (\hat{e}, \cdots, \hat{e} ).$$

\noindent
{\it Proof of Theorem 4.2.} Let $[u]$ be one of generators with $u=z f_j^{(n)} + (1-f_j^{(n)})$. First, we compute $\psi([u]) = \langle \varepsilon\sharp\tau, u\rangle=\frac{1}{2\pi i}\varepsilon\sharp\tau(u^{-1}, u)$. As $u^{-1}=u^* = \overline{z}f+(1-f)$ where $f=f_j^{(n)}$, we have
$$ \frac{1}{2\pi i}\varepsilon\sharp\tau(u^{-1}, u) = \frac{1}{2\pi i}\varepsilon\sharp\tau(\overline{z}f+(1-f),z f+(1-f)) = \frac{1}{2\pi i}\varepsilon(\overline{z}, z)\tau(f^2)=\tau(f).$$
In order to compute $\langle \omega, \delta_1([u])\rangle$ we apply the Falk Index Formula. 

Let $F\in C(T_0)$ be an extension of $f = f_j^{(n)}$. Set $M=LF+(1-F), N=FL^*+ (1-F) \in A$. Then $\pi(M)= u, \pi(N) = u^*$. For $a, b \in B(H)$, if $a - b\in \K(H)$ we denote by $a\equiv b$. We have $[L, F] \equiv 0, F(1-F) \equiv 0$. From this
$$1-MN \equiv F(1-LL^*), 1-NM \equiv 0.$$
Hence $(1-MN)^2\equiv F(1-LL^*)^2, (1-NM)^2\equiv 0$. Therefore
$$\pi(\hat{e})= \begin{pmatrix}
f(1-MM^*)^2 & \ast\ast \\
\ast & 0 \\
\end{pmatrix}.$$
Notice that $1-MM^*$ is a rank one projection on $\ell^2(\Z_+)$. This implies that
$$\langle \omega, \delta_1([u])\rangle=\omega(\hat{e}) = \textup{Tr}(1-MM^*)\tau(f) = \tau(f).$$
\begin{flushright}
$\square$
\end{flushright}

\noindent
{\it Proof of Corollary 4.4.} By Theorem 4.2 it is sufficient to compute $\langle \;\varepsilon\sharp\tau, [\hat{\sigma}(Q_+)] \;\rangle$. For any given $x\in K$ the function $\hat{\sigma}(Q_+)(x, \cdot)$ is invertible in $C(S^1)$. We compute $\langle \varepsilon, [\hat{\sigma}(Q_+)(x, \cdot)] \rangle$.

Set $\alpha=\frac{\sqrt{1+a(x)}\,}{\sqrt{1-a(x)}}\overline{b}(x), \beta=\frac{\sqrt{1-a(x)}}{\sqrt{1+a(x)}\,}b(x)$, and set
$g(z) = \hat{\sigma}(Q_+)(x,z) = \alpha z-\beta \overline{z}.$

Suppose that $a(x) > 0$. Then $|\frac{\beta}{\,\alpha\,}| = \frac{1-a(x)}{\,1+a(x)\,}<1$, and
\begin{equation*}
\begin{split}
g^{-1} & = \frac{1}{\alpha z-\beta \overline{z}} = \frac{\overline{z}}{\alpha- \beta \overline{z}^2} 
 = \frac{ \overline{z}}{\;\alpha\;}\cdot \frac{1}{1-(\beta/\alpha) \overline{z}^2} \\
& =  \frac{ \overline{z}}{\,\alpha\,}\sum_{n=0}^\infty \bigl(\tfrac{\beta}{\,\alpha\,}\overline{z}^2\bigr)^n.
\end{split}
\end{equation*}
Hence
$$g^{-1} dg = \frac{ \overline{z}}{\,\alpha\,}\sum_{n=0}^\infty \bigl(\tfrac{\beta}{\,\alpha\,}\overline{z}^2\bigr)^n (\alpha dz- \beta d\overline{z}).$$
From this we have
$$\int_{S^1}g^{-1} dg = \int_{S^1} \tfrac{ \overline{z}}{\,\alpha\,}\sum_{n=0}^\infty \bigl(\tfrac{\beta}{\,\alpha\,}\overline{z}^2\bigr)^n \cdot\alpha dz = \int_{S^1} \tfrac{ \overline{z}}{\,\alpha\,}\cdot \alpha dz= 2\pi i.$$

When $a(x) < 0$, we have $|\frac{\alpha}{\,\beta\,}| < 1$. By arguments similar to the above, we have
$$\int_{S^1}g^{-1} dg = \int_{S^1} \tfrac{z}{\,\beta\,}\sum_{n=0}^\infty \bigl(\tfrac{\alpha}{\,\beta\,}z^2\bigr)^n \cdot (-\beta) dz = \int_{S^1} \tfrac{z}{\,\beta\,} (-\beta d\overline{z})= -2\pi i.$$
Therefore
$$\langle \varepsilon, [\hat{\sigma}(Q_+)(x, \cdot)] \rangle = \left \{ \begin{array}{cl}
 1 & \mbox{if } a(x) > 0 \\
 -1 & \mbox{if } a(x) < 0\\
 \end{array}
 \right. .$$
 From this the conclusion follows.
 \begin{flushright}
 $\square$
  \end{flushright}
 \section{a generalization}
 
 As in the classical setting we can introduce further parameters. Fix $(p, q) \in S^2 \subset \R\times\C$. Then 
 $$\Gamma = \begin{pmatrix}
p  & \overline{q}L^* \\
qL & E-p \\
\end{pmatrix}$$
is a symmetry. We use the coin operator $C$ defined above. Then the chirality operator is
$$Q_+ = \frac{\overline{b}}{\sqrt{1-a}}\,qL\sqrt{1+a} - \sqrt{1-a}\,\overline{q}L^*\frac{b}{\sqrt{1+a}} +\frac{\overline{b}}{\sqrt{1-a}}E\frac{b}{\sqrt{1+a}} -2p|b|.$$
For each $x\in K, z\in S^1$ we have
\begin{equation*}
\begin{split}
\hat{\sigma}(Q_+)(x,z) &= \frac{q\overline{b(x)}(1+a(x))}{|b(x)|}z- \frac{\overline{q}b(x)(1-a(x))}{\sqrt{1-a(x)^2}}\overline{z}-2p|b(x)|\\
& = q(1+a(x))w -\overline{q}(1-a(x))\overline{w}-2p|b(x)|\\
& = g(x, w),
\end{split}
\end{equation*}
where $b(x)=|b(x)|e^{i\theta}$ and $w=e^{-i\theta}z$.

For a fixed $x\in K$ consider a linear equation $q(1+a(x))w -\overline{q}(1-a(x))\overline{w}-2p|b(x)| =0$ in $w$. The unique solution is
$$w_0 =\frac{\overline{q}p|b(x)|}{\,\,a(x)|q|^2\,\,}.$$
It is easy to show that $|w_0| =1$ if and only if $p^2=a(x)^2$. Thus

\begin{prop}
The symbol $\hat{\sigma}(Q_+)$ is invertible in $C(S^1)\otimes C(K)$ if and only if $|p|\neq |a(x)|$ for all $x\in K$.
\end{prop}

For a fixed $x\in K$ we want to compute 
\begin{equation*}
\begin{split}
\langle \varepsilon, [\hat{\sigma}(Q_+)(x,\cdot)]\rangle & = \frac{1}{\,2\pi i\,}\int_{S^1} \frac{1}{\hat{\sigma}(Q_+)(x,z)}d\hat{\sigma}(Q)(x,z)\\
& =  \frac{1}{\,2\pi i\,}\int_{S^1} \frac{1}{g(x, w)}d g(x,w).
\end{split}
\end{equation*}

We have
$$g^{-1}dg = \frac{w^2q(1+a(x))+\overline{q}(1-a(x))}{w\bigl(q(1+a(x))w^2-2p|b(x)|w-\overline{q}(1-a(x)\bigr)\bigr)} dw.$$

The method we used above does not seem plausible to be applied to our new situation. So we apply a classical method, {\it i.e.} the residue formula. The integrand is a meromorphic function with three poles $0, \alpha, \beta$ where
$$\alpha = \frac{|b(x)|(p+1)}{q(1+a(x))}, \;\beta=\frac{|b(x)|(p-1)}{q(1+a(x))}.$$
By straightforward computations we get the residues at $0, \alpha, \beta$ as follows:
$$Res(0) =-1, Res(\alpha) =1, Res(\beta)=1.$$
It is easy to see that
(1)\; $|\alpha| = 1$ if and only if $a(x)= p$, 
(2) \,$|\alpha| < 1 $ if and only if $a(x) > p$,
and (3) \,$|\alpha| > 1 $ if and only if $a(x) <  p$.
Similarly
(1)\; $|\beta| = 1$ if and only if $a(x)= -p$, 
(2) \,$|\beta| < 1 $ if and only if $a(x) > -p$,
and (3) \,$|\beta| > 1 $ if and only if $a(x) <  -p$.

Under the assumption $p >0$, if $a(x) > p$, then $\alpha, \beta \in \D$. It follows that
$$I = \frac{1}{\,2\pi i\,}\int_{S^1} \frac{1}{g(x, w)}d g(x,w)= Res(0)+Res(\alpha)+Res(\beta)=1.$$
If $-p < a(x) < p$, then $\alpha\notin \D, \beta \in \D$ . Hence
$$I = Res(0) + Res(\beta)=0.$$
If $a(x) < -p$, we have $\alpha, \beta \notin \D$. Then
$$I = Res(0) =-1.$$
Let $p < 0$. If $a(x) > -p$, then $\alpha, \beta \in \D$. It follows that
$$I = \frac{1}{\,2\pi i\,}\int_{S^1} \frac{1}{g(x, w)}d g(x,w)= Res(0)+Res(\alpha)+Res(\beta)=1.$$
If $p < a(x) < -p$, then $\alpha\in \D, \beta \notin \D$ . So
$$I = Res(0) + Res(\alpha)=0.$$
If $a(x) < p$, then $\alpha, \beta \notin \D$. Consequently
$$I = Res(0) =-1.$$

Set $K_+ = \{ x\in K : a(x) >|p| \}, K_0=\{ x\in K : -|p| < a(x) < |p| \},$ and $K_- = \{ x\in K : a(x) <-|p| \}$. Denote by $\chi_+,  \chi_0,  \chi_-$ the characteristic functions of $K_+, K_0, K_-$, respectively.

Summarizing the arguments above:
 \begin{thm}
We have  $\textup{s-ind} \,Q_+ = \tau(\chi_+) - \tau(\chi_-)$.
\end{thm}

\begin{remark}
In the arguments above we did not discuss the domains of cyclic cocycles. However, it is straightforward to see that the domains are large enough to detect $K$-theory.
\end{remark}

\end{document}